\documentclass[a4paper,11pt]{article}
\usepackage{epsfig,makeidx,amsmath,amsfonts,latexsym, subfigure}
\newtheorem{theorem}{Theorem}[section]
\newtheorem{prop}{Proposition}[section]
\newtheorem{lemma}{Lemma}[section]

\def\0{{\bf 0}}
\def\P{{\bf P}}
\def\E{{\bf E}}

\begin{document}

\title{Random intersection graphs with tunable degree
distribution and clustering}

\author{Maria Deijfen\thanks{Stockholm University.
E-mail: mia@math.su.se} \and Willemien Kets\thanks{Santa Fe
Institute and Tilburg University. E-mail:
willemien.kets@santafe.edu}}

\date{August 6, 2008}

\maketitle

\thispagestyle{empty}

\begin{abstract}
\noindent A random intersection graph is constructed by assigning independently to each vertex a subset of a given set and
drawing an edge between two vertices if and only if their
respective subsets intersect. In this paper a model is
developed in which each vertex is given a random weight, and
vertices with larger weights are more likely to be assigned
large subsets. The distribution of the degree of a given vertex
is characterized and is shown to depend on the weight of the
vertex. In particular, if the weight distribution is a power
law, the degree distribution will be so as well. Furthermore,
an asymptotic expression for the clustering in the graph is
derived. By tuning the parameters of the model, it is possible
to generate a graph with arbitrary clustering, expected degree
and -- in the power law case -- tail exponent.

\vspace{0.5cm}

\noindent \emph{Keywords:} Random intersection graphs, degree
distribution, power law distribution, clustering, social
networks.

\vspace{0.5cm}

\noindent AMS 2000 Subject Classification: 05C80, 91D30.
\end{abstract}

\section{Introduction}

During the last decade there has been a large interest in the
study of large complex networks; see e.g. Dorogovtsev and
Mendes (2003) and Newman et al.\ (2006) and the references
therein. Due to the rapid increase in computer power, it has
become possible to investigate various types of real networks
such as social contact structures, telephone networks, power
grids, the Internet and the World Wide Web. The empirical
observations reveal that many of these networks have similar
properties. For instance, they typically have power law degree
sequences, that is, the fraction of vertices with degree $k$ is
proportional to $k^{-\tau}$ for some exponent $\tau>1$.
Furthermore, many networks are highly clustered, meaning
roughly that there is a large number of triangles and other
short cycles. In a social network, this is explained by the
fact that two people who have a common friend often meet and
become friends, creating a triangle in the network. A related
explanation is that human populations are typically divided
into various subgroups -- working places, schools, associations
etc -- which gives rise to high clustering in the social
network, since members of a given group typically know each
other; see Palla et al. (2005) for some empirical observations.

Real-life networks are generally very large, implying that it
is a time-consuming task to collect data to delineate their
structure in detail. This makes it desirable to develop models
that capture essential features of the real networks. A natural
candidate to model a network is a random graph, and, to fit
with the empirical observations, such a graph should have a
heavy-tailed degree distribution and considerable clustering.
We will quantify the clustering in a random graph by the
conditional probability that three given vertices constitute a
triangle, given that two of the three possible links between
them exist. Other (empirical) definitions occur in the
literature -- see e.g.\ Newman (2003) -- but they all capture
essentially the same thing.

Obviously, the classical Erd\H{o}s-R\'{e}nyi graph will not do
a good job as a network model, since the degrees are
asymptotically Poisson distributed. Moreover, existing models
for generating graphs with a given degree distribution -- see
e.g.\ Molloy and Reed (1995, 1998) -- typically have zero
clustering in the limit. In this paper, we propose a model,
based on the so-called random intersection graph, where both
the degree distribution and the clustering can be controlled.
More precisely, the model makes it possible to obtain arbitrary
prescribed values for the clustering and to control the mean
and the tail behavior of the degree distribution.

\subsection{Description of the model}

The random intersection graph was introduced in Singer (1995)
and Karo\'{n}ski et\ al.\ (1999), and has been further studied
and generalized in Fill et\ al. (2000), Godehardt and Jaworski
(2002), Stark (2004) and Jaworksi et\ al.\ (2006). Newman
(2003) and Newman and Park (2003) discuss a similar model. In
its simplest form the model is defined as follows.
\begin{enumerate}

\item Let $\mathcal{V}=\{1,\ldots,n\}$ be a set of $n$
    vertices and $\mathcal{A}$ a set of $m$ elements. For
    $p\in[0,1]$, construct a bipartite graph $B(n,m,p)$
    with vertex sets $\mathcal{V}$ and $\mathcal{A}$ by
    including each one of the $nm$ possible edges between
    vertices from $\mathcal{V}$ and elements from
    $\mathcal{A}$ independently with probability $p$.

\item The random intersection graph $G(n,m,p)$ with vertex
    set $\mathcal{V}$ is obtained by connecting two
    distinct vertices $i,j \in \mathcal{V}$ if and only if
    there is an element $a\in \mathcal{A}$ such that both
    $i$ and $j$ are adjacent to $a$ in $B(n,m,p)$.
\end{enumerate}
When the vertices in $\mathcal{V}$ are thought of as
individuals and the elements of $\mathcal{A}$ as social groups,
this gives rise to a model for a social network in which two
individuals are joined by an edge if they share at least one
group. In the following, we frequently borrow the terminology
from the field of social networks and refer to the vertices as
individuals and the elements of $\mathcal{A}$ as groups, with
the understanding that the model is of course much more
general.

To get an interesting structure, the number of groups $m$ is
typically set to $m=\lfloor n^\alpha\rfloor$ for some $\alpha >
0$; see Karo\'{n}ski et\ al.\ (1999). We will assume this form
for $m$ in the following. Let $D_i$ be the degree of vertex $i
\in \mathcal{V}$ in $G(n,m,p)$. The probability that two
individuals do not share a group in $B(n,m,p)$ is $(1-p^2)^m$.
It follows that the edge probability in $G(n,m,p)$ is
$1-(1-p^2)^m$ and hence the expected degree is
\begin{eqnarray*}
\E[D_i] & = & (n-1)(1-(1-p^2)^m)\\
& = & (n-1)\left(mp^2+O(m^2p^4)\right).
\end{eqnarray*}
To keep the expected degree bounded as $n\to\infty$, we let
$p=\gamma n^{-(1+\alpha)/2}$ for some constant $\gamma>0$. We
then have that $\E[D_i]\to \gamma^2$.

Stark (2004; Theorem 2) shows that in a random intersection
graph with the above choice of $p$, the distribution of the
degree of a given vertex converges to a point mass at 0, a
compound Poisson distribution or a Poisson distribution
depending on whether $\alpha<1$, $\alpha=1$ or $\alpha>1$. This
means that the current model cannot account for the power law
degree distributions typically observed in real networks.

In the above formulation of the model, the number of groups
that a given individual belongs to is binomially distributed
with parameters $m$ and $p$. A generalization of the model,
allowing for an arbitrary group distribution, is described in
Godehardt and Jaworski (2002). The degree of a given vertex in
such a graph is analyzed in Jaworski et al.\ (2006), where
conditions on the group distribution are specified under which
the degree is asymptotically Poisson distributed.

In the current paper, we are interested in obtaining graphs
where non-Poissonian degree distributions can be identified. To
this end, we propose a generalization of the original random
intersection graph where the edge probability $p$ is random and
depends on weights associated with the vertices. Other work in
this spirit include for instance Chung and Lu (2002:1,2), Yao
et\ al.\ (2005), Britton et\ al.\ (2006), Bollob\'{a}s et\ al.\
(2007) and Deijfen et\ al.\ (2007). The model is defined as
follows:

\begin{enumerate}
\item Let $n$ be a positive integer, and define
    $m=\lfloor\beta n^\alpha\rfloor$ with $\alpha,\beta>0$.
    As before, take $\mathcal{V}=\{1,\ldots,n\}$ to be a
    set of $n$ vertices and $\mathcal{A}$ a set of $m$
    elements. Also, let $\{W_i\}$ be an i.i.d.\ sequence of
    positive random variables with distribution $F$, where
    $F$ is assumed to have mean 1 if the mean is finite.
    Finally, for some constant $\gamma>0$, set
\begin{equation}\label{p_i}
p_i=\gamma W_i n^{-(1+\alpha)/2}\wedge 1.
\end{equation}
Now construct a bipartite graph $B(n,m,F)$ with vertex sets
$\mathcal{V}$ and $\mathcal{A}$ by adding edges to the
elements of $\mathcal{A}$ for each vertex $i \in
\mathcal{V}$ independently with probability $p_i$.

\item The random intersection graph $G(n,m,F)$ is obtained
    as before by drawing an edge between two distinct
    vertices $i,j \in \mathcal{V}$ if and only if they have
    a common adjacent vertex $a \in \mathcal{A}$ in
    $B(n,m,F)$.
\end{enumerate}

In the social network setting, the weights can be interpreted
as a measure of the social activity of the individuals. Indeed,
vertices with large weights are more likely to join many groups
and thereby acquire many social contacts. There are several
other examples of real networks where the success of a vertex
(measured by its degree) depends on some specific feature of
the vertex; see e.g.\ Palla et al.\ (2005) for an example in
the context of protein interaction networks. Furthermore, an
advantage of the model is that it has an explicit and
straightforward construction which, as we will see, makes it
possible to exactly characterize the degree distribution and
the clustering in the resulting graph.

\subsection{Results}

Our results concern the degree distribution and the clustering
in the graph $G(n,m,F)$ as $n\to\infty$. More precisely, we
will take the parameters $\alpha$, $\beta$, $\gamma$ and the
weight distribution $F$ to be fixed (independent of $n$) and
then analyze the degree of a given vertex and the clustering in
the graph as $n\to\infty$. It turns out that the behavior of
these quantities will be different in the three regimes
$\alpha<1$, $\alpha=1$ and $\alpha>1$ respectively. The
interesting case is $\alpha=1$, in the sense that this is when
both the degree distribution and the clustering can be
controlled. The cases $\alpha<1$ and $\alpha>1$ are included
for completeness.

As for the degree, we begin by observing that, if $F$ has
finite mean, then the asymptotic mean degree of vertex $i$,
conditional on $W_i$, is given by $\beta\gamma^2W_i$ for all
values of $\alpha$.

\begin{prop}\label{prop:vv} Let $D_i$ be the degree of vertex $i\in\mathcal{V}$ in a random intersection graph $G(n,m,F)$ with $m=\lfloor \beta n^\alpha \rfloor$ and $p_i$ as in~\eqref{p_i}. If $F$ has finite mean, then, for all values of $\alpha>0$, we have that $\E[D_i|W_i]\to\beta\gamma^2W_i$ as $n\to\infty$.
\end{prop}

\noindent \textbf{Proof.} We prove the claim for vertex $i=1$.
Define
\[
W_j' = W_j\cdot \textbf{1}_{\{W_j\leq n^{1/4}\}} \quad\textrm{and}\quad W_j'' = W_j\cdot\textbf{1}_{\{W_j> n^{1/4}\}}
\]
and let $D'$ and $D''$ denote the degree of vertex 1 when
$\{W_j\}_{j\neq 1}$ are replaced by $\{W_j'\}$ and $\{W_j''\}$
respectively, that is, $D'$ is the number of neighbors of 1
with weight smaller than or equal to $n^{1/4}$ and $D''$ is the
number of neighbors with weight larger than $n^{1/4}$. Write
$p_j'$ and $p_j''$ for the analog of (\ref{p_i}) based on the
truncated weights.

Now, conditional on the weights, the probability that there is
an edge between $1$ and $j$ is $1-(1-p_1p_j)^m$. To see that
$\E[D'']\to 0$ as $n\to\infty$, we observe that
\[
1-(1-p_ip_j'')^m\leq mp_1p_j''=\beta\gamma W_1n^{(\alpha-1)/2}p_j''.
\]
Summing the expectation of the right-hand side over $j\neq 1$,
keeping $W_1$ fixed, gives (recall the truncation at 1 in
(\ref{p_i}))
\[
\E[D''] \leq \beta\gamma n^{(1+\alpha)/2}\E[p_k'']\leq \beta\gamma\left(\gamma\E[W_k'']+n^{(1+\alpha)/2}\P\big(\gamma W_k\geq n^{(1+\alpha)/2}\big)\right),
\]
where both terms on the right hand side converge to 0 as
$n\to\infty$ since $F$ has finite mean. As for $D'$, we have
\[
1-(1-p_1p_j')^m=\beta\gamma^2W_1W_j'n^{-1}+O(W_1^2(W_j')^2n^{-2}).
\]
The sum over $j\neq 1$ of the expectation of the first term
equals $\beta\gamma^2 W_1\E[W_k']$, where
$\E[W_k']\to\E[W_k]=1$ (since $F$ has finite mean) and the sum
of the expectation of the second term converges to 0 (since
$(W_j')^2\leq n^{1/2}$). Since $D_0=D'+D''$, this proves the proposition.\hfill$\Box$\medskip

The following theorem, which is a generalization of Theorem 2
in Stark (2004), gives a full characterization of the degree
distribution for different values of $\alpha$.

\begin{theorem}\label{th:deg_distr} Consider the degree $D_i$ of vertex
$i \in \mathcal{V}$ in a random intersection graph $G(n,m,F)$
with $m=\lfloor \beta n^\alpha \rfloor$ and $p_i$ as
in~\eqref{p_i}, and assume that $F$ has finite mean.

\begin{itemize}
\item[\rm{(a)}] If $\alpha<1$, then  $D_i$ converges in
    distribution to a point mass at 0 as $n\to\infty$.

\item[\rm{(b)}] If $\alpha=1$, then $D_i$ converges in
    distribution to a sum of a Poisson($\beta\gamma W_i$)
    distributed number of Poisson($\gamma$) variables,
    where all variables are independent.

\item[\rm{(c)}] If $\alpha>1$, then $D_i$ is asymptotically
    Poisson($\beta\gamma^2W_i$) distributed.
\end{itemize}
\end{theorem}

To understand Theorem \ref{th:deg_distr}, note that the
expected number of groups that individual $i$ belongs to is
roughly $\beta\gamma W_in^{(\alpha-1)/2}$. If $\alpha<1$ and
$W_i$ has finite mean, this converges to 0 in
probability, so that the degree distribution converges to a
point mass at 0, as stated in (a) (the group size however goes
to infinity, explaining why the expected degree is still
positive in the limit). For $\alpha=1$, the number of groups
that individual $i$ is a member of is Poisson($\beta\gamma
W_i$) distributed as $n\to\infty$, and the number of other
individuals in each of these groups is approximately
Poisson($\gamma$) distributed, which explains (b). Finally, for
$\alpha>1$, individual $i$ belongs to infinitely many groups as
$n\to\infty$. This means that the edges indicators will be
asymptotically independent, giving rise to the Poisson
distribution specified in (c).

Moving on to the clustering, write $E_{ij}$ for the event that
individuals $i,j\in\mathcal{V}$ have a common group in the
bipartite graph $B(n,m,F)$ -- that is, $E_{ij}$ is equivalent
to the event that there is an edge between vertices $i$ and $j$
in $G(n,m,F)$ -- and let $\bar{\P}_n$ be the probability
measure of $B(n,m,F)$ conditional on the weights
$\{W_1,\ldots,W_n\}$. For distinct vertices
$i,j,k\in\mathcal{V}$, define

\begin{equation}\label{c_n}
\bar{c}_{i,j,k}^{(n)}=\bar{\P}_n\left(E_{ij}|E_{ik},E_{jk}\right),
\end{equation}

\noindent that is, $\bar{c}_{i,j,k}^{(n)}$ is the edge
probability between $i$ and $j$ in $G(n,m,F)$ given that they
are both connected to $k$, conditional on the weights. To
quantify the asymptotic clustering in the graph we will use
$$
c(G):=\lim_{n\to\infty}\E\big[\bar{c}_{i,j,k}^{(n)}\big],
$$
where the expectation is taken over the weights, that
is, $c(G)$ is the limiting probability that three given
vertices constitute a triangle conditional on that two of the
three possible edges between them exist (the vertices are
indistinguishable, so indeed $c(G)$ does not depend on the
particular choice of $i,j$ and $k$). This should be closely
related to the limiting quotient of the number of triangles and
the number of triples with at least two edges present, which is
one of the empirical measures of clustering that occur in the
literature; see e.g.\ Newman (2003). Establishing this
connection rigorously however requires additional
arguments.\medskip

The asymptotic behavior of $\bar{c}_{i,j,k}^{(n)}$ is specified
in the following theorem. By bounded convergence it follows
that $c(G)$ is obtained as the mean of the
in-probability-limits.

\begin{theorem}\label{th:clustering} Let $\{i,j,k\}$ be three distinct vertices in a random intersection graph $G(n,m,F)$ with $m=\lfloor\beta n^\alpha\rfloor$ and $p_i$
as in (\ref{p_i}). If $F$ has finite mean, then
\begin{itemize}
\item[\rm{(a)}] $\bar{c}_{i,j,k}^{(n)}\to 1$ in probability
    for $\alpha<1$;

\item[\rm{(b)}] $\bar{c}_{i,j,k}^{(n)}\to(1+\beta\gamma
    W_k)^{-1}$ in probability for $\alpha=1$;

\item[\rm{(c)}] $\bar{c}_{i,j,k}^{(n)}\to 0$ in probability
    for $\alpha>1$.
\end{itemize}
\end{theorem}

To understand Theorem~\ref{th:clustering}, assume that $i$ and
$k$ share a group and that $j$ and $k$ share a group. The
probability that $i$ and $j$ also have a common group then
depends on the number of groups that the common neighbor $k$
belongs to. Indeed, the fewer groups $k$ belongs to, the more
likely it is that $i$ and $j$ in fact share the same group with
$k$. Recall that the expected number of groups that $k$ belongs
to is roughly $\beta\gamma W_kn^{(\alpha-1)/2}$. If $\alpha>1$,
this goes to 0 as $n\to\infty$. Since it is then very unlikely
that $k$ belongs to more than one group when $n$ is large, two
given edges $\{i,k\}$ and $\{j,k\}$ are most likely generated
by the same group, meaning that $i$ and $j$ are connected as
well. On the other hand, if $\alpha>1$, the number of groups
that $k$ belongs to is asymptotically infinite. Hence, that $i$
and $j$ each belong to one of these groups, does not
automatically make it likely that they actually belong to the
same group. If $\alpha=1$, individual $k$ belongs to
$\beta\gamma W_k$ groups on average, explaining the expression
in part (b) of the theorem.

From Theorem \ref{th:clustering} it follows that, to get a
nontrivial tunable clustering, we should choose $\alpha=1$.
Indeed, then we have $c(G)=\E[(1+\beta\gamma W_k)^{-1}]$ and,
for a given weight distribution $F$ (with finite mean), $c(G)$
can be varied between 0 and 1 by adjusting the parameters
$\beta$ and $\gamma$. Furthermore, when $\alpha=1$, the degree
distribution for a given vertex is asymptotically compound
Poisson with the weight of the vertex as one of the parameters
-- see Theorem \ref{th:deg_distr} (b) -- and it is not hard to
see that, if $F$ is a power law with exponent $\tau$, then the
degree distribution will be so as well. Since the mean of $F$
is set to 1, the expected asymptotic degree is $\beta\gamma^2$
by Proposition \ref{prop:vv}. Taken together, this means that,
when $\alpha=1$, we can obtain a graph with a given value of
the clustering and a power law degree distribution with
prescribed exponent and prescribed mean by first choosing $F$
to be a power law with the desired exponent and then tuning the
parameters $\beta$ and $\gamma$ to get the correct values of
the clustering and the expected degree.

The rest of the paper is organized as follows. In Sections
\ref{proof:deg distr} and \ref{proof:clustering},
Theorem~\ref{th:deg_distr} and Theorem~\ref{th:clustering} are
proved, respectively. The clustering is analyzed for the
important example of a power law weight distribution in
Section~\ref{power law}. Finally, Section~\ref{conclusions}
provides an outline of possible future work.

\section{The degree distribution} \label{proof:deg distr}

We begin by proving Theorem \ref{th:deg_distr}.\medskip

\noindent \textbf{Proof of Theorem \ref{th:deg_distr}.} We
prove the theorem for vertex $i=1$. Write $D_1=D$, and denote
by $N$ the number of groups that individual 1 belongs to.
Conditional on $W_1$, the variable $N$ is binomially
distributed with parameters $m$ and $p_1$ and thus
\[
\bar{\P}_n(N=0)=(1-p_1)^m\geq 1-mp_1\geq 1-\beta\gamma^2 W_1^{(\alpha-1)/2}.
\]
For $\alpha<1$, the expectation of the last term converges to 0
as $n\to\infty$, and it follows from bounded convergence that
$\P(N=0)=\E[\bar{\P}_n(N=0)]\to 1$. This proves (a), since
clearly $D=0$ if individual 1 is not a member of any group.

To prove (b) and (c), first recall the definition of the
weights $\{W_i'\}$ and $\{W_i''\}$ -- truncated from above and
below respectively at $n^{1/4}$ -- and the corresponding degree
variables $D'$ and $D''$ from the proof of Proposition
\ref{prop:vv}. We have already showed (in proving Proposition
\ref{prop:vv}) that $\E[D'']\to 0$, which implies that $D''$
converges to 0 in probability (indeed $P(D''>0)\leq \E[D'']$).
Hence it suffices to show that the generating function of $D'$
converges to the generating function of the claimed limiting
distribution. To this end, we condition on the weight $W_1$,
which is thus assumed to be fixed in what follows, and let
$X_i'$ ($i=2,\ldots,n$) denote the number of common groups of
individual $1$ and individual $i$ when the truncated weights
$W_i'$ are used for $i\neq 1$. Since two individuals are
connected if and only if they have at least one group in
common, we can write $D'=\sum_{i=2}^n\mathbf{1}_{\{X_i'\geq
1\}}$. Furthermore, conditional on $N$ and $\{W_i'\}_{i\geq
2}$, the random variables $X_i'$, $i =2, \ldots, n$, are
independent and binomially distributed with parameters $N$ and
$p_i'=\gamma W_i'n^{-(1+\alpha)/2}$. Hence, with
$\bar{\bar{\P}}_n$ denoting the probability measure of the
bipartite graph $B(n,m,F)$ conditional on both $\{W_i'\}_{i\geq
2}$ and $N$, the generating function of $D'$ can be written as
\[
\E[t^{D'}] = \E\left[\prod_{i=2}^n\E\left[t^{\mathbf{1}\{X_i'\geq
1\}}\big|\{W_i'\},N\right]\right] = \E\left[\prod_{i=2}^n\left(1+(t-1)\bar{\bar{\P}}_n(X_i'\geq
1)\right)\right]
\]
where $t\in[0,1]$. Using the Taylor expansion
$\log(1+x)=x+O(x^2)$ and the fact that
\[
\bar{\bar{\P}}_n(X_i'\geq 1)=1-(1-p_i')^N=Np_i'+O(N^2(p_i')^2),
\]
we get that
\begin{eqnarray}\label{gen_fct_prod}
\prod_{i=2}^n\left(1+(t-1)\bar{\bar{\P}}_n(X_i'\geq 1)\right) & = &
e^{(t-1)N\sum p_i'+O\big(N^2\sum (p_i')^2\big)}.
\end{eqnarray}

\noindent Defining
\[
R_n:=\prod_{i=2}^n\left(1+(t-1)\bar{\bar{\P}}_n(X_i'\geq
1)\right)-e^{(t-1)N\sum p_i'}
\]
we therefore have that
\[
R_n=e^{(t-1)N\sum p_i}\Big(e^{O(N^2\sum (p_i')^2)}-1\Big).
\]
Since the product in (\ref{gen_fct_prod}) is the conditional
expectation of $t^{D'}$ with $t\in[0,1]$, it takes values
between 0 and 1 and, since $e^{(t-1)N\sum p_i}\in (0,1]$, it
follows that $R_n\in[-1,1]$. Furthermore, recalling that
$W_i'\leq n^{1/4}$, we have for $\alpha\geq 1$ that
\[
N^2\sum_{i=2}^n(p_i')^2=N^2\gamma^2n^{-(1+\alpha)}\sum_{i=2}^n(W_i')^2\leq N^2\gamma^2n^{-1/2},
\]
implying that $R_n\to 0$ in probability and thus, by bounded
convergence, $\E[R_n]\to 0$. Hence we are done if we show that
\begin{itemize}
\item[{\rm{(i)}}] $\E\left[e^{(t-1)N\sum p_i}\right]\to
    e^{\beta\gamma W_1(e^{\gamma(t-1)}-1)}$ if $\alpha=1$;

\item[{\rm{(ii)}}] $\E\left[e^{(t-1)N\sum p_i}\right]\to
    e^{\beta\gamma^2 W_1(t-1)}$ if $\alpha>1$,
\end{itemize}
where the limits are recognized as the generating functions for
the desired compound Poisson and Poisson distribution in part
(b) and (c) of the theorem, respectively. To this end, note
that the expectation with respect to $N$ of $e^{(t-1)N\sum
p_i'}$ is given by the generating function for $N$ evaluated at
the point $e^{(t-1)\sum p_i'}$. Since $N$ is binomially
distributed with parameters $m$ and $p_1$, we have that
\begin{equation}\label{(i,ii)gen_fct}
\E\left[e^{(t-1)N\sum p_i'}\right]=\E\left[\left(1+p_1\left(e^{(t-1)\sum p_i'}-1\right)
\right)^m\right].
\end{equation}
For $\alpha=1$, we have $m=\lfloor \beta n\rfloor$ and
$p_i'=\gamma W_i'n^{-1}$. Recalling that $\E[W_i']\to
\E[W_i]=1$, it follows that $\sum p_i'\to \gamma$ almost
surely. Hence,
\[
\left(1+p_1\left(e^{(t-1)\sum p_i'}-1\right)\right)^{\lfloor \beta
n\rfloor}\to e^{\beta\gamma W_1(e^{\gamma(t-1)}-1)}\textrm{ a.s.\ as
}n\to\infty,
\]
and it follows from bounded convergence that the expectation
converges to the same limit, proving (i).

For $\alpha>1$, define $\tilde{p}_i'=n^{(\alpha-1)/2}p_i'$.
With $m=\lfloor\beta n^\alpha\rfloor$ and $p_1=\gamma
W_1n^{-(1+\alpha)/2}\wedge 1$, we get after some rewriting,
that
\[
\left(1+p_1\left(e^{(t-1)\sum p_i'}-1\right)\right)^m=
\left(1+\frac{\gamma W_1(t-1)\sum\tilde{p}_i'}{n^\alpha}\cdot\frac{e^{(t-1)n^{(1-\alpha)/2}\sum
\tilde{p}_i'}-1}{(t-1)n^{(1-\alpha)/2}\sum\tilde{p}_i'}\right)^{\lfloor\beta
n^\alpha\rfloor}.
\]
By the law of large numbers, $\sum\tilde{p}_i'\to\gamma$ almost
surely, and, since $(e^x-1)/x\to 1$ as $x\to 0$, it follows
that the right hand side above converges to
$e^{\beta\gamma^2W_1(t-1)}$ almost surely as $n\to\infty$. By
(\ref{(i,ii)gen_fct}) and bounded convergence, this proves
(ii).\hfill$\Box$

\section{Clustering} \label{proof:clustering}

In this section, we prove Theorem \ref{th:clustering}. First
recall that $E_{ij}$ denotes the event that the individuals
$i,j\in\mathcal{V}$ share at least one group. It will be
convenient to extend this notation. To this end, for $i,j,k \in
\mathcal{V}$, denote by $E_{ijk}$ the event that there is at
least one group to which all three individuals $i$, $j$ and $k$
belong, and write $E_{ij,ik,jk}$ for the event that there are
at least three \emph{distinct} groups to which $i$ and $j$, $i$
and $k$, and $j$ and $k$ respectively belong. Similarly, the
event that there are two distinct groups to which individuals
$i$ and $k$, and $j$ and $k$ respectively belong is denoted by
$E_{ik,jk}$. The proof of Theorem \ref{th:clustering} relies on
the following lemma.

\begin{lemma} \label{lemma:clustering} Consider a random intersection
graph $G(n,m,F)$ with $m=\lfloor\beta n^\alpha\rfloor$ and
$p_i$ defined as in (\ref{p_i}). For any three distinct
vertices $i,j,k\in\mathcal{V}$, we have that
\begin{itemize}
\item[\rm{(a)}] $\bar{\P}_n(E_{ijk})= \frac{\beta
    \gamma^3W_iW_jW_k}{n^{(3+\alpha)/2}}+O\left(\frac{W_i^2W_j^2W_k^2}{n^{3+\alpha}}\right)$;

\item[\rm{(b)}] $\bar{\P}_n(E_{ij,ik,jk})=
    \frac{\beta^3\gamma^6W_i^2W_j^2W_k^2}{n^3}+
    O\left(\frac{W_i^3W_j^3W_k^3}{n^4}\right)$;

\item[\rm{(c)}]
    $\bar{\P}_n(E_{ik,jk})=\frac{\beta^2\gamma^4W_iW_jW_k^2}{n^2}+
    O\left(\frac{W_i^2W_j^2W_k^3}{n^3}\right)$;

\item[\rm{(d)}]
    $\bar{\P}_n(E_{ijk}E_{ik,jk})=O\left(\frac{W_i^2W_j^2W_k^2}{n^{(5+\alpha)/2}}\right)$.
\end{itemize}
\end{lemma}

\noindent \textbf{Proof.} As for (a), the probability that
three given individuals $i$, $j$ and $k$ do not share any group
at all is $(1-p_ip_jp_k)^m$. Using the definitions of $m$ and
the edge probabilities $\{p_i\}$, it follows that
\[
\bar{\P}_n(E_{ijk}) = 1-(1-p_ip_jp_k)^m =  \frac{\beta\gamma^3W_iW_jW_k}{n^{(3+\alpha)/2}}+
O\left(\frac{W_i^2W_j^2W_k^2}{n^{3+\alpha}}\right).
\]
To prove (b), note that the probability that there is exactly
one group to which both $i$ and $j$ belong is
$mp_ip_j(1-p_ip_j)^{m-1}=mp_ip_j+O(m^2p_i^2p_j^2)$. Given that
$i$ and $j$ share one group, the probability that $i$ and $k$
share exactly one of the \emph{other} $m-1$ groups is
$(m-1)p_ip_k(1-p_ip_k)^{m-2}=mp_ip_k+O(m^2p_i^2p_k^2)$.
Finally, the conditional probability that there is a third
group to which both $j$ and $k$ belong given that the pairs
$i,j$ and $i,k$ share one group each is
$1-(1-p_jp_k)^{m-2}=mp_jp_k+O(m^2p_j^2p_k^2)$. Combining these
estimates, and noting that scenarios in which $i$ and $j$ or
$i$ and $k$ share more than one group have negligible
probability in comparison, we get that
\begin{eqnarray*}
\bar{\P}_n(E_{ij,ik,jk}) & = & m^3p_i^2p_j^2p_k^2+O\left(m^4p_i^2p_j^2p_k^2(p_ip_j+p_ip_k+p_jp_k)\right)\\
& = & \frac{\beta^3\gamma^6W_i^2W_j^2W_k^2}{n^3}+
O\left(\frac{W_i^3W_j^3W_k^3}{n^4}\right).
\end{eqnarray*}
Part (c) is derived analogously.

As for (d), note that the event $E_{ijk}E_{ik,jk}$ occurs when
there is at least one group that is shared by all three
vertices $i$, $j$ and $k$ and a second group shared by either
$i$ and $k$ or $j$ and $k$. Denote by $r$ the probability that
individual $k$ and at least one of the individuals $i$ and $j$
belong to a fixed group. Then $r=p_k(p_i+p_j-p_ip_j)$, and,
conditional on that there is exactly one group to which all
three individuals $i$, $j$ and $k$ belong (the probability of
this is $mp_ip_jp_k(1-p_ip_jp_k)^{m-1}=O(mp_ip_jp_k)$), the
probability that there is at least one other group that is
shared either by $i$ and $k$ or by $j$ and $k$ is
$1-(1-r)^{m-1}=O(mr)$. It follows that
\[
\bar{\P}_n(E_{ijk}E_{ik,jk}) = O(m^2p_ip_jp_kr) =  O\left(\frac{W_i^2W_j^2W_k^2}{n^{(5+\alpha)/2}}\right).
\]
\hfill$\Box$\medskip

Using Lemma \ref{lemma:clustering}, it is not hard to prove
Theorem \ref{th:clustering}.\medskip

\noindent \textbf{Proof of Theorem \ref{th:clustering}.} Recall
the definition (\ref{c_n}) of $\bar{c}_{i,j,k}^{(n)}$ and note
that
\[
\bar{\P}_n(E_{ij}|E_{ik}E_{jk})=\frac{\bar{\P}_n(E_{ijk}\cup
E_{ij,ik,jk})}{\bar{\P}_n(E_{ijk}\cup E_{ik,jk})}.
\]
As for (a), applying the estimates of Lemma
\ref{lemma:clustering} and merging the error terms yields
\begin{eqnarray}
\bar{\P}_n(E_{ij}|E_{ik}E_{jk}) & \geq &\nonumber
\frac{\bar{\P}_n(E_{ijk})}{\bar{\P}_n(E_{ijk})+\bar{\P}_n(E_{ik,jk})}\\
& = & \frac{1+O(W_iW_jW_kn^{-(3+\alpha)/2})}{1+W_k[\beta\gamma
n^{(\alpha-1)/2}+O(W_iW_jW_kn^{-(3-\alpha)/2})]}.\label{lower_bound}
\end{eqnarray}
By Markov's inequality and the fact that $W_i$, $W_j$ and $W_k$
are independent and have finite mean, it follows that
$W_iW_jW_kn^{-(3-\alpha)/2}$ goes to 0 in probability when
$\alpha<1$. Similarly, $W_iW_jW_kn^{-(3+\alpha)/2}\to 0$ in
probability. Hence, the quotient in (\ref{lower_bound})
converges to 1 in probability for $\alpha<1$, as claimed.

To prove part (b), note that, for $\alpha=1$, the lower bound
(\ref{lower_bound}) for $\bar{c}_{i,j,k}^{(n)}$ converges in
probability to $(1+\beta\gamma W_k)^{-1}$. To obtain an upper
bound, we apply Lemma \ref{lemma:clustering} with $\alpha=1$ to
get that
\begin{eqnarray}
\bar{\P}_n(E_{ij}|E_{ik}E_{jk}) & \leq &
\frac{\bar{\P}_n(E_{ijk})+\bar{\P}_n(E_{ij,ik,jk})}{\bar{\P}_n(E_{ijk})+\bar{\P}_n(E_{ik,jk})-
\bar{\P}_n(E_{ijk}E_{ik,jk})}\label{cl_bd}\\
& = &
\frac{1+O(W_iW_jW_kn^{-1})}{1+W_k[\beta\gamma+O(W_iW_jW_kn^{-1})]}.\nonumber
\end{eqnarray}
Here $W_iW_jW_kn^{-1}$ converges to 0 in probability by
Markov's inequaliy, and (b) follows.

As for (c), combining the bound in (\ref{cl_bd}) with the
estimates in Lemma \ref{lemma:clustering} yields
\[
\bar{\P}_n(E_{ij}|E_{ik}E_{jk}) \leq \frac{n^{(1-\alpha)/2}+O(W_iW_jW_kn^{-1})} {n^{(1-\alpha)/2}+W_k[\beta\gamma+O(W_iW_jW_kn^{-1})]}.
\]
Since $W_iW_jW_kn^{-1}\to 0$ in probability, this bound
converges to 0 in probability for $\alpha>1$, as
desired.\hfill$\Box$

\section{Clustering for a power law weight distribution} \label{power law}

When $\alpha=1$, the clustering is given by
$c(G)=\E[(1+\beta\gamma W_k)^{-1}]$. Here we investigate this
expression in more detail for the important case that $F$ is a
power law. More precisely, we take $F$ to be a Pareto
distribution with density
\[
f(x)=\frac{(\lambda-2)^{\lambda-1}}{(\lambda-1)^{\lambda-2}}x^{-\lambda}\quad\textrm{for
}x\geq \frac{\lambda-2}{\lambda-1}.
\]
When $\lambda>2$, this distribution has mean 1, as desired. The
asymptotic clustering $c(G)$ is given by the integral
$$
\frac{(\lambda-2)^{\lambda-1}}{(\lambda-1)^{\lambda-2}}\int_{\frac{\lambda-2}{\lambda-1}}^\infty
(1+\beta\gamma x)^{-1}x^{-\lambda}dx.
$$
Defining $u:= (\lambda - 2)/(x \cdot (\lambda -1))$, we obtain
\begin{eqnarray*}
c(G) &=& \frac{1}{\beta\gamma} \frac{(\lambda-1)^2}{(\lambda-2)}
\int_0^1 u^{\lambda -1} \left(1 + \frac{u}{\beta
\gamma}\left( \frac{\lambda -1}{\lambda - 2}\right) \right)^{-1} du\\
&=:& \frac{1}{\beta\gamma\lambda}
\frac{(\lambda-1)^2}{(\lambda-2)}  \mbox{ }
{}_2F_1\left(1,\lambda; 1 + \lambda; -\frac{1}{\beta \gamma}
\left(\frac{\lambda -1}{\lambda -2}\right)\right),
\end{eqnarray*}
where ${}_2F_1$ is the hypergeometric function. For $\beta \gamma \geq (\lambda -1)/(\lambda
- 2)$, a series expansion of the integrand yields that
\begin{eqnarray*}
c(G) &= &\frac{1}{\beta\gamma} \frac{(\lambda-1)^2}{(\lambda-2)}
\sum_{k = 0}^\infty \left( - \frac{1}{\beta \gamma}\left(
\frac{\lambda -1}{\lambda - 2}\right) \right)^k \frac{1}{k +
\lambda}\\
&=:& \frac{1}{\beta\gamma} \frac{(\lambda-1)^2}{(\lambda-2)}
\Phi\left(- \frac{1}{\beta \gamma}\left( \frac{\lambda -1}{\lambda
- 2}\right),  1, \lambda \right),
\end{eqnarray*}
where $\Phi$ is the Lerch transcedent. Furthermore, when
$\lambda$ is an integer, we get
\begin{eqnarray*}
c(G) &=& \frac{(\lambda-2)^{\lambda-1}}{(\lambda-1)^{\lambda-2}}
\left[ (-\beta \gamma)^{\lambda - 1} \ln\left( 1+\frac{\lambda
-1}{\beta \gamma (\lambda -2)}\right) + \sum_{\ell = 1}^{\lambda -
1}
 \frac{(-\beta \gamma)^{\lambda - 1 -
\ell}}{\ell }  \left(\frac{\lambda - 1}{\lambda-2} \right)^\ell
\right].
\end{eqnarray*}

\begin{figure}
\centering \mbox{\subfigure[The clustering as a function of $\lambda$ for
different values of $\beta \gamma$: $\beta\gamma =1$ (---),
$\beta\gamma =5$ ($- - -$), $\beta\gamma =10$ ($- \cdot -$).]{\epsfig{file=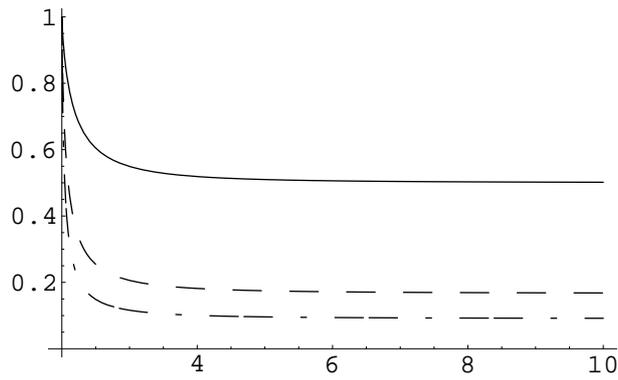,height=5cm}}}\par
\mbox{\subfigure[The clustering as a function of $\beta \gamma$ for different
values of $\lambda$: $\lambda =2.1$ (---), $\lambda =2.5$ ($- -
-$), $\lambda =4$ ($- \cdot -$).]{\epsfig{file=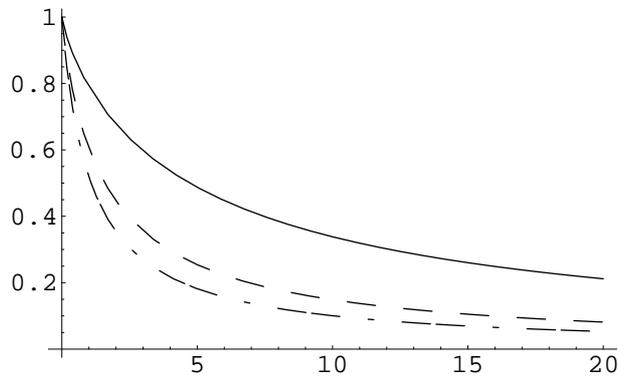,height=5cm}}}\caption{Clustering for a power law distribution.}\label{graph_clust}
\end{figure}

Figure~\ref{graph_clust} (a) and (b) show how the clustering
depends on $\lambda$ and $\beta\gamma$ respectively. For any $c
\in (0,1)$ and a given tail exponent $\lambda$, we can find a
value of $\beta\gamma$ such that the clustering is equal to
$c$. Combining this with a condition on $\beta\gamma^2$,
induced by fixing the mean degree in the graph, the parameters
$\beta$ and $\gamma$ can be specified.

\section{Future work} \label{conclusions}

Apart from the degree distribution and the clustering, an
important feature of real networks is that there is typically
significant correlation for the degrees of neighboring nodes,
that is, either high (low) degree vertices tend to be connected
to other vertices with high (low) degree (positive
correlation), or high (low) degree vertices tend to be
connected to low (high) degree vertices (negative correlation).
A next step is thus to quantify the degree correlations in the
current model. The fact that individuals share groups should
indeed induce positive degree correlation, which agrees with
empirical observations from social networks; see Newman (2003)
and Newman and Park (2003).

Also other features of the model are worth investigating. For
instance, many real networks are ``small worlds'', meaning
roughly that the distances between vertices remain small also
in very large networks. It would be interesting to study the
relation between the distances between vertices, the degree
distribution and the clustering in the current model.

Finally, dynamic processes behave differently on clustered
networks as compared to more tree-like networks. Most work to
date has focused on the latter class. In Britton et\ al.\
(2008) however, epidemics on random intersection graphs without
random weights are studied and it is investigated how the
epidemic spread is affected by the clustering in the graph. It
would be interesting to extend this work to incorporate weights
on the vertices, allowing to tune also the (tail of the) degree
distribution and study its impact on the epidemic process.\medskip

\noindent \textbf{Acknowledgement.} We thank Remco van der Hofstad and
Wouter Kager for valuable suggestions that have improved the manuscript.

\section*{References}

\noindent Bollob\'{a}s, B., Janson, S. and Riordan, O. (2007): The
phase transition in inhomogeneous random graphs, \emph{Random
Structures $\&$ Algorithms} \textbf{31}, 3-122.\medskip

\noindent Britton, T., Deijfen, M. and Martin-L\"{o}f, A.
(2006): Generating simple random graphs with prescribed degree
distribution, \emph{Journal of Statistical Physics} \textbf{124},
1377-1397.\medskip

\noindent Britton, T., Deijfen, M., Lager\aa s, A. and
Lindholm, M. (2008): Epidemics on random graphs with tunable
clustering, \emph{Journal of Applied Probability} \textbf{45}, 743-756.\medskip

\noindent Chung, F. and Lu. L (2002:1): Connected components in
random graphs with given degree sequences, \emph{Annals of
Combinatorics} \textbf{6}, 125-145.\medskip

\noindent Chung, F. and Lu. L (2002:2): The average distances
in random graphs with given expected degrees, \emph{Proceedings
of the National Academy of Sciences} \textbf{99},
15879-15882.\medskip

\noindent Deijfen, M., vd Esker, H., vd Hofstad, R. and
Hooghiemstra, G. (2007): A preferential attachment model with
random initial degrees, \emph{Arkiv f\"{o}r Matematik}, to
appear.\medskip

\noindent Dorogovtsev, S. and Mendes, J. (2003):
\emph{Evolution of Networks, from Biological Nets to the
Internet and WWW}, Oxford University Press.\medskip

\noindent Fill, J., Scheinerman, E. and Singer-Cohen, K.
(2000): Random intersection graphs when $m=\omega(n)$: an
equivalence theorem relating the evolution of the $G(n,m,p)$
and $G(n,p)$ models, \emph{Random Structures $\&$ Algorithms}
\textbf{16}, 156-176.\medskip

\noindent Godehardt, E. and Jaworski, J. (2002): Two models of
random intersection graphs for classification, in
\emph{Exploratory data analysis in empirical research}, eds.\
Schwaiger M. and Opitz, O., Springer, 67-81.\medskip

\noindent Jaworski, J., Karo\'{n}ski, M. and Stark, D. (2006):
The degree of a typical vertex in generalized random
intersection graph models, \emph{Discrete Mathematics}
\textbf{306}, 2152-2165.\medskip

\noindent Karo\'{n}ski, M., Scheinerman, E. and Singer-Cohen,
K. (1999): On random intersection graphs: the subgraphs
problem, \emph{Combinatorics, Probability $\&$ Computing}
\textbf{8}, 131-159.\medskip

\noindent Molloy, M. and Reed, B. (1995): A critical point for
random graphs with a given degree sequence, \emph{Random
Structures $\&$ Algorithms} \textbf{6}, 161179.\medskip

\noindent Molloy, M. and Reed, B. (1998): The size of the giant
component of a random graph with a given degree sequence,
\emph{Combinatorics, Probability $\&$ Computing} \textbf{7},
295-305.\medskip

\noindent Newman, M. E. J., Strogatz, S. H., and Watts, D. J.
(2002): Random graphs with arbitrary degree distributions and
their applications, \emph{Physical Review E} \textbf{64},
026118.\medskip

\noindent Newman, M. E. J. (2003): Properties of highly
clustered networks, \emph{Physical Review E} \textbf{68},
026121.\medskip

\noindent Newman, M. E. J. and Park J. (2003): Why social
networks are different from other types of networks,
\emph{Physical Review E} \textbf{68}, 036122.\medskip

\noindent Palla, G., and Der{\'e}nyi, I., Farkas, I. and
Vicsek, T. (2005): Uncovering the overlapping community
structure of complex networks in nature and society,
\emph{Nature} \textbf{435}, 814 - 818. \medskip

\noindent Singer, K. (1995): \emph{Random intersection graphs},
PhD thesis, Johns Hopkins University.\medskip

\noindent Stark, D. (2004): The vertex degree distribution of
random intersection graphs, \emph{Random Structures $\&$
Algorithms} \textbf{24}, 249-258.\medskip

\noindent Yao, X., Zhang, C., Chen, J. and Li, Y. (2005): On
the scale-free intersection graphs, \emph{Lecture notes in
computer science} \textbf{3481}, 1217-1224.

\vfill\eject

\end{document}